\documentclass[a4paper]{amsart}
\usepackage[]{fontenc}
\usepackage[latin1]{inputenc}
\usepackage{verbatim}
\usepackage{amssymb}

\makeatletter
  \theoremstyle{plain}
  \newtheorem{lem}{Lemma}
  \theoremstyle{plain}
  \newtheorem{thm}{Theorem}

\usepackage{hyperref}
\usepackage{mathrsfs}

\usepackage[english]{babel}
\makeatother
\begin{document}
\selectlanguage{english}
\newcommand{\set}[2]{\left\{  #1\, |\, #2\right\}  }

\newcommand{\cyc}[1]{\mathbb{Q}\left[\zeta_{#1}\right]}

\newcommand{\ZN}[1]{\left(\mathbb{Z}/#1\mathbb{Z}\right)^{\times}}

\newcommand{\Mod}[3]{#1\equiv#2\, \left(\mathrm{mod}\, \, #3\right)}

\newcommand{\sm}[4]{\left(\begin{smallmatrix}#1  &  #2\cr\cr#3  &  #4\end{smallmatrix}\right)}

\newcommand{\aut}[1]{\mathrm{Aut}\!\left(#1\right)}

\newcommand{\End}[1]{\mathrm{End}\!\left(#1\right)}

\newcommand{\gl}[2]{\mathrm{GL}_{#1}\left(#2\right)}

\newcommand{\tr}[2]{\textrm{Tr}_{_{#1}}\left(#2\right)}

\newcommand{\FA}{\mathbb{Z}\oplus\mathbb{Z}}

\newcommand{\FB}{\mathbf{H}}

\newcommand{\FC}{\mathrm{SL}_{2}\!\left(\mathbb{Z}\right)}

\newcommand{\FD}{\ell_{i}}

\newcommand{\FE}[1]{\mathscr{L}\!\left(#1\right)}

\newcommand{\FF}[1]{\mathscr{C}\!\left(#1,R\right)}

\newcommand{\FG}[1]{#1^{*}}

\newcommand{\FH}[2]{}

\newcommand{\FI}[2]{}

\newcommand{\FJ}[2]{}

\title{The orbifold transform and its applications}

\author{P. Bantay}

\email{bantay@general.elte.hu}

\keywords{finitely generated groups, permutation orbifolds, wreath products}

\thanks{Work supported by grants OTKA T047041, T043582, the J{\' a}nos Bolyai
Research Scholarship of the Hungarian Academy of Sciences and EC Marie
Curie MRTN-CT-2004-512194.}

\curraddr{Institute for Theoretical Physics, E{\" o}tv{\" o}s Lor{\' a}nd University, Budapest}

\begin{abstract}
We discuss the notion of the orbifold transform, and illustrate it
on simple examples. The basic properties of the transform are presented,
including transitivity and the exponential formula for symmetric products.
The connection with the theory of permutation orbifolds is addressed,
and the general results illustrated on the example of torus partition
functions.
\end{abstract}
\maketitle

\section{Introduction}

In the last few decades, Conformal Field Theory (CFT) \cite{BPZ,DiFrancesco-Mathieu-Senechal}
and the closely related String Theory \cite{GSW,Polch} have had many
fruitful interactions with different branches of mathematics, ranging
from algebraic topology through differential geometry to the theory
of modular forms, not to mention the intimately related theories of
Vertex Operator Algebras \cite{FLM1,Kac} and Modular Tensor Categories\cite{Turaev,Bakalov-Kirillov}.

Orbifold constructions, i.e. the gauging of discrete internal symmetries,
have played an important role in Conformal Field Theory for quite
some time \cite{Dixon_orbifolds1,Dixon_orbifolds2,Dixon_orbifoldCFT}.
Among these, the theory of permutation orbifolds \cite{Klemm-Schmidt,Borisov-Halpern-Schweigert,permorbi1,permorbi2}
is a most interesting chapter: not only does it provide a general
procedure for constructing new Conformal Field Theories from known
ones, with pretty good control over the structure of the resulting
theory, but -- through the so-called Orbifold Covariance Principle
\cite{OCP} -- one has a very effective tool for the study of deeper
aspects of CFT, which led ultimately to a proof of the  Congruence
Subgroup conjecture for Rational CFT \cite{CMP} (a related proof
in the context of nets of subfactors has been recently provided in
\cite{Xu}). Symmetric product orbifolds, i.e. permutation orbifolds
of the full symmetric group play a basic role in the description of
second quantized strings \cite{elliptic_genera,Dijkgraaf_disctors,symprod,LuninMathur3}.

Many important aspects of permutation orbifolds can be understood
to a large extent with the help of a general group theoretic construct
that we term the orbifold transform. In special instances one gets
back well known classical concepts, like the cycle indicator polynomial
of finite permutation groups, while in other special cases the transform
describes correlation and partition functions of permutation orbifolds.
The aim of the present paper is describe this general notion and to
discuss its properties, as well as its applications to CFT.

In the next section, we'll introduce the concept of the orbifold transform,
and illustrate it in some simple cases. Section 3 is devoted to the
statement and proof of a most important property of the orbifold transform,
namely transitivity, which describes the result of successive applications
of the transform. Section 4 is concerned with the proof and applications
of a general combinatorial identity satisfied by the orbifold transform,
which plays a fundamental role in the theory of symmetric products.
In Section 5 we sketch the connection of the orbifold transform with
the theory of permutation orbifolds. Finally, we conclude by indicating
possible extensions, which might prove usefull in describing such
concepts as discrete torsion. A short appendix treats some combinatorial
results used in the text.

\section{The orbifold transform}

To start with, let's recall the following basic facts about finitely
generated groups and their finite index subgroups \cite{Magnus-Karras-Solitar,Robinson}:

\begin{itemize}
\item a finite index subgroup of a finitely generated group is itself finitely
generated (this follows from the Reidemeister-Schreier theorem);
\item a finite index subgroup of a finite index subgroup is again of finite
index (by Lagrange's theorem), and the intersection of two finite
index subgroups is again of finite index (by a theorem of Poincar�;
\item there are only finitely many different homomorphisms from a finitely
generated group into a finite group (since a homomorphism is determined
by the images of the generators, and there are only finitely many
possibilities for them);
\item a finitely generated group has only finitely many conjugacy classes
of subgroups of a given finite index (because each conjugacy class
corresponds to a transitive permutation action of degree equal to
the index, and these form a finite set by the above), and each such
conjugacy class contains finitely many different subgroups (because
the normalizer of a finite index subgroup is obviously also of finite
index).
\end{itemize}
\medskip{}
\medskip{}

Let $G$ be a finitely generated group, and let $\FE{G}$ denote the
set of finite index subgroups of $G$. Let $R$ be a commutative $\mathbb{Q}$-algebra,
i.e. a commutative ring with identity that contains an isomorphic
copy of the field of rational numbers. By a class function with values
in $R$ we shall mean a mapping $Z:\FE{G}\rightarrow R$ which is
constant on conjugacy classes of subgroups, i.e. such that \begin{equation}
Z\left(g^{-1}Hg\right)=Z\left(H\right)\:\label{eq:classfun}\end{equation}
holds for all $g\in G$ and $H\in\FE{G}$. We shall denote by $\FF{G}$
the set of such class functions with values in $R$.

Given a permutation group $\Omega<S_{X}$ acting on the finite set
$X$, we can associate to any $Z\in\FF{G}$ a new map \begin{equation}
\begin{aligned}\begin{aligned}Z\wr\Omega:\FE{G}\rightarrow\end{aligned}
 & \: R\\
H\mapsto & \:\frac{1}{\left|\Omega\right|}\sum_{\phi:H\rightarrow\Omega}\:\prod_{\xi\in\mathcal{O}\left(\phi\right)}Z\left(H_{\xi}\right)\:.\end{aligned}
\label{eq:orbitransdef}\end{equation}

Let's see the different ingredients entering this definition! First,
the summation runs over all homomorphisms $\phi:H\rightarrow\Omega$
mapping the subgroup $H$ into $\Omega$: since $H$ is finitely generated
and $\Omega$ is finite, the sum is finite. For each homomorphism
$\phi:H\rightarrow\Omega$, the image $\phi\left(H\right)$ -- being
a subgroup of $\Omega$, hence of $S_{X}$ -- is a permutation group
acting on the set $X$, and we denote by $\mathcal{O}\left(\phi\right)$
the set of its orbits: there are only finitely many of them, since
$X$ is finite. Finally, for a given orbit $\xi\in\mathcal{O}\left(\phi\right)$
we denote by $H_{\xi}$ any of its point stabilizers, more precisely
 \begin{equation}
H_{\xi}=\set{g\in H}{\xi^{*}\:\mathrm{is}\:\mathrm{fixed}\:\mathrm{by}\:\phi\left(g\right)}\:\label{eq:stabdef}\end{equation}
for some representative $\xi^{*}\in\xi$ chosen at will: $H_{\xi}$
is a finite index subgroup, its index in $H$ being equal to the length
$\left|\xi\right|$ of the orbit $\xi$. Note that, although the subgroup
$H_{\xi}$ does depend on the choice of the representative $\xi^{*}$,
the value $Z\left(H_{\xi}\right)$ doesn't, since stabilizers of points
on the same orbit are conjugate subgroups and $Z$ is a class function.

In summary, every term in Eq.\eqref{eq:orbitransdef} makes sense:
$Z\left(H_{\xi}\right)$ exists and is independent of the chosen representative
$\xi^{*}\in\xi$; the sum and the product are both finite; and we
can divide by the order $\left|\Omega\right|$ of the permutation
group $\Omega$, since $R$ is a $\mathbb{Q}$-algebra. Consequently,
Eq.\eqref{eq:orbitransdef} gives a well defined map $Z\wr\Omega:\FE{G}\rightarrow R$,
which we call the orbifold transform of the class function $Z$.

A basic property of the orbifold transform is that $Z\wr\Omega$ is
itself a class function. To see this, let's consider a conjugate $K=H^{g}$
of the subgroup $H\in\FE{G}$, and let's determine the value of $Z\wr\Omega$
on $K$: according to Eq.\eqref{eq:orbitransdef}, this is given by

\[
\frac{1}{\left|\Omega\right|}\sum_{\phi:K\rightarrow\Omega}\:\prod_{\xi\in\mathcal{O}\left(\phi\right)}Z\left(K_{\xi}\right)\:.\]
But to each $\phi:H\rightarrow\Omega$ we can associate a homomorphism
$\phi^{g}:K\rightarrow\Omega$ via the rule $\phi^{g}\left(h^{g}\right)=\phi\left(h\right)$,
and this correspondence is clearly one-to-one. Moreover, the image
of $\phi^{g}$ equals that of $\phi$, so they have the same orbits,
and the stabilizers of the orbit $\xi\in\mathcal{O}\left(\phi\right)=\mathcal{O}\left(\phi^{g}\right)$
are conjugate subgroups: $K_{\xi}=g^{-1}H_{\xi}g$. Since $Z$ is
a class function, we conclude that $Z\wr\Omega$ agrees on $H$ and
$K$, i.e. it is a class function too, as claimed.

In summary, for each permutation group $\Omega$ of finite degree
the orbifold transform provides a map \begin{eqnarray*}
\FF{G} & \rightarrow & \FF{G}\\
Z & \mapsto & Z\wr\Omega\qquad.\end{eqnarray*}
To get some familiarity with this map, let's see how it looks like
in simple cases!

Of course, the simplest case is when $G$ is trivial, i.e.  consists
of the identity element solely. There is just one subgroup of $G=\left\{ 1\right\} $,
namely $G$ itself, so a class function is nothing but an element
of $R$. There is just one homomorphism $\phi:G\rightarrow\Omega$
for any permutation group $\Omega$, namely the one that takes the
identity of $G$ to the identity of $\Omega$, and each point of $X$
forms an orbit in itself, whose stabilizer is the whole of $G$. Consequently,
in this particular case the orbifold transform takes the form \begin{equation}
Z\wr\Omega=\frac{1}{\left|\Omega\right|}Z^{d}\:,\label{eq:trivigroup}\end{equation}
where $d$ denotes the degree of $\Omega$, i.e. the cardinality of
$X$.

One can analyze in a similar fashion the case of any finite $G$:
the orbifold transform gives then -- for each permutation group $\Omega$
of finite degree -- a polynomial map $R\left[t_{1},\ldots,t_{n}\right]\rightarrow R\left[t_{1},\ldots,t_{n}\right]$,
where $n$ denotes the number of conjugacy classes of subgroups of
$G$, whose explicit expression depends on both $G$ and $\Omega$.

Let's turn to a less trivial case, when $G=\mathbb{Z}$ is infinite
cyclic! Since $G$ is Abelian, each conjugacy class consists of one
subgroup; moreover, for each positive integer $n$ there is exactly
one subgroup of index $n$, namely $G_{n}=n\mathbb{Z}$, and these
are all infinite cyclic. This means that a class function $Z\in\FF{G}$
may be viewed as an infinite sequence $z_{1},z_{2},\ldots$ of elements
of $R$, where $z_{n}=Z\left(G_{n}\right)$. Since $G_{n}$ is generated
by one element, a homomorphism $\phi:G_{n}\rightarrow\Omega$ is specified
by giving the image of the generator, which could be an arbitrary
element $x$ of $\Omega$: the image of $\phi$ is then the cyclic
subgroup of $\Omega$ generated by $x$, and the orbits of the image
are nothing but the orbits of $x$. Finally, the stabilizer of a given
orbit $\xi$ -- which is independent of the choice of the representative
$\xi^{*}\in\xi$, since the subgroup $G_{n}$ is Abelian -- is $\set{g\in G_{n}}{g^{\left|\xi\right|}=1}$,
which is nothing else but the subgroup $G_{n\left|\xi\right|}$ of
$G$. Summarizing all this, we have\begin{equation}
Z\wr\Omega:G_{n}\mapsto\frac{1}{\left|\Omega\right|}\sum_{x\in\Omega}\:\prod_{\xi\in\mathcal{O}\left(x\right)}Z\left(G_{n\left|\xi\right|}\right)\:.\label{eq:Ztransform}\end{equation}

We see that the orbifold transform maps the infinite sequence $z_{1},z_{2},\ldots$
into another infinite sequence $z_{1}^{\Omega},z_{2}^{\Omega},\ldots\in R$,
whose elements are finite polynomials in the $z_{i}$. Actually, this
map has a classical interpretation: to see this, recall from enumerative
combinatorics \cite{stanley} the notion of the cycle indicator of
the permutation group $\Omega$, which is the multivariate polynomial\begin{equation}
P_{\Omega}\left(t_{1},\ldots,t_{d}\right)=\frac{1}{\left|\Omega\right|}\sum_{x\in\Omega}\:\prod_{\xi\in\mathcal{O}\left(x\right)}t_{\left|\xi\right|}\;,\label{eq:cyleindicator}\end{equation}
where $d$ denotes the degree of $\Omega$ and the $t_{i}$ are indeterminates.
With the help of the cycle indicator, we can recast Eq.\eqref{eq:Ztransform}
for the orbifold transform into\begin{equation}
z_{n}^{\Omega}=P_{\Omega}\left(z_{n},z_{2n},\ldots,z_{dn}\right)\:.\label{eq:Ztransform2}\end{equation}

Finally, let's consider a most interesting case for applications in
the theory of permutation orbifolds, when $G=\mathbb{Z}\oplus\mathbb{Z}$
is free Abelian of rank two (this is the fundamental group of a two-dimensional
torus, explaining its relevance to Conformal Field Theory). Again,
such a $G$ is Abelian, so each conjugacy class contains just one
subgroup. Finite index subgroups of $\mathbb{Z}\oplus\mathbb{Z}$
are all isomorphic to $\mathbb{Z}\oplus\mathbb{Z}$, and they are
in one-to-one correspondence with 2-by-2 integer matrices in Hermite
normal form (HNF) \cite{Cohen1}. Recall, that such a matrix has the
form \begin{equation}
H=\left(\begin{array}{cc}
\mu & \kappa\\
0 & \lambda\end{array}\right)\:,\label{eq:HNF}\end{equation}
where $\mu$ and $\lambda$ are positive integers, while $\kappa$
is a nonnegative integer less than $\lambda$, i.e. $0\leq\kappa<\lambda$.
If $G=\mathbb{Z}\oplus\mathbb{Z}$ is generated by $a$ and $b$ (note
that $a$ and $b$ commute), then the finite index subgroup corresponding
to $H$ is generated by $a^{\lambda}$ and $a^{\kappa}b^{\mu}$, and
its index equals the determinant of $H$. From now on we shall freely
identify the matrix $H$ with the corresponding subgroup, and view
class functions $Z\in\FF{G}$ as defined on the set of matrices in
HNF.

A homomorphism $\phi:\FA\rightarrow\Omega$ is specified by a pair
$\left(x,y\right)$ of commuting elements of $\Omega$ (the images
of the generators): the image of such a homomorphism is the subgroup
of $\Omega$ generated by $x$ and $y$, and we shall denote by $\mathcal{O}\left(x,y\right)$
the orbits of this subgroup. The stabilizer of each orbit $\xi\in\mathcal{O}\left(x,y\right)$
-- which is once again independent of the choice of the representative
$\xi^{*}\in\xi$, since the group  is Abelian -- is a finite index
subgroup of $\FA$, so it corresponds to a matrix \[
H_{\xi}=\left(\begin{array}{cc}
\mu_{\xi} & \kappa_{\xi}\\
0 & \lambda_{\xi}\end{array}\right)\]
in HNF: here, $\lambda_{\xi}$ denotes the common length of all $x$
orbits contained in $\xi$ and $\mu_{\xi}$ denotes their number,
while $\kappa_{\xi}$ characterizes the 'skewness' of the orbit $\xi$.
Taking all this into account, the orbifold transform reads \begin{equation}
Z\wr\Omega:H\mapsto\frac{1}{\left|\Omega\right|}\sum_{{{x,y\in\Omega\atop xy=yx}}}\:\prod_{\xi\in\mathcal{O}\left(x,y\right)}Z\left(H_{\xi}H\right)\:.\label{eq:torustransform}\end{equation}
When suitably interpreted, the above formula gives the torus partition
function of the permutation orbifold with twist group $\Omega$ \cite{permorbi1}.

Hopefully, the above examples were able to give an impression of the
general notion of the orbifold transform and its relation to some
classical notions. But it should be stressed that this construct works
for any finitely generated group, leading to some genuinely new structures
in general.

\section{Transitivity of the orbifold transform}

By far the most important property of the orbifold transform is transitivity,
which describes the result of successive applications of the transform,
and is closely related to the corresponding property of permutation
orbifolds \cite{permorbi1,permorbi2}. Since the formulation of this
property, as well as its proof, relies strongly on the theory of wreath
products, let's begin by recalling some facts about the latter \cite{dixon-mortimer,Kerber}.

Consider two permutation groups $\Omega_{1}<S_{X}$ and $\Omega_{2}<S_{Y}$
 acting on the finite sets $X$ and $Y$. To any $\lambda\in\Omega_{1}^{Y}$
-- mapping $Y$ to $\Omega_{1}$ -- and a permutation $\omega\in\Omega_{2}$
one can associate a permutation $\lambda\wr\omega$ of $X\times Y$
via the rule \begin{equation}
\lambda\wr\omega:\left(x,y\right)\mapsto\left(\lambda\left(y\right)x,\omega y\right)\;.\label{eq:wreathaction}\end{equation}
The important observation is that the product of two permutations
of this type is again of this type, to wit \begin{equation}
\left(\lambda_{1}\wr\omega_{1}\right)\left(\lambda_{2}\wr\omega_{2}\right)=\left(\lambda_{1}^{\omega_{2}}\lambda_{2}\right)\wr\left(\omega_{1}\omega_{2}\right)\:,\label{eq:wreathmult}\end{equation}
where $\lambda_{1}^{\omega_{2}}\lambda_{2}:Y\rightarrow\Omega_{1}$
is given by \[
\lambda_{1}^{\omega_{2}}\lambda_{2}:y\mapsto\lambda_{1}\left(\omega_{2}y\right)\lambda_{2}\left(y\right)\:.\]
In short, the permutations in $\Omega_{1}\wr\Omega_{2}=\set{\lambda\wr\omega}{\lambda\in\Omega_{1}^{Y},\:\omega\in\Omega_{2}}$
form a group that acts on the Cartesian product $X\times Y$,  called
the wreath product of $\Omega_{1}$ with $\Omega_{2}$. Clearly, the
degree of $\Omega_{1}\wr\Omega_{2}$ is the product of the degrees
of its factors: \begin{equation}
\deg\left(\Omega_{1}\wr\Omega_{2}\right)=\deg\left(\Omega_{1}\right)\deg\left(\Omega_{2}\right)\:,\label{eq:wreathdegree}\end{equation}
while the order of the wreath product is given by \begin{equation}
\left|\Omega_{1}\wr\Omega_{2}\right|=\left|\Omega_{1}\right|^{\deg\left(\Omega_{2}\right)}\left|\Omega_{2}\right|\:.\label{eq:wreathorder}\end{equation}
The orbits of the wreath product are easy to describe: each orbit
of $\Omega_{1}\wr\Omega_{2}$ is of the form $\xi\times\eta$, where
$\xi$ is an orbit of $\Omega_{1}$ on $X$  while $\eta$ is an orbit
of $\Omega_{2}$ on $Y$.

Wreath products are associative but not commutative, in the sense
that the permutation groups $\left(\Omega_{1}\wr\Omega_{2}\right)\wr\Omega_{3}$
and $\Omega_{1}\wr\left(\Omega_{2}\wr\Omega_{3}\right)$ are always
equivalent, but $\Omega_{1}\wr\Omega_{2}$ and $\Omega_{2}\wr\Omega_{1}$
fail to be equivalent in general.

Consider a homomorphism $\phi:G\rightarrow\Omega_{1}\wr\Omega_{2}$
from an arbitrary group $G$ into the wreath product $\Omega_{1}\wr\Omega_{2}$.
Such a $\phi$ assigns to each element $g\in G$ a permutation $\phi\left(g\right)=\lambda\left(g\right)\wr\omega\left(g\right)$.
This means that $\phi$ can be described by the pair $\left(\lambda,\omega\right)$,
where $\lambda$ maps $G$ into $\Omega_{1}^{Y}$, while $\omega$
maps $G$ into $\Omega_{2}$. Taking into account Eq.\eqref{eq:wreathmult},
one sees that actually $\omega:G\rightarrow\Omega_{2}$ is a homomorphism,
while $\lambda:G\rightarrow\Omega_{1}^{Y}$ is a ($\omega$-)crossed
homomorphism, i.e. a map that satisfies  \begin{equation}
\lambda\left(gh\right)=\lambda\left(g\right)^{\omega\left(h\right)}\lambda\left(h\right)\:.\label{eq:crossedhomo}\end{equation}

To formulate the next result that classifies crossed homomorphisms,
let's recall that for an orbit $\xi\in\mathcal{O}\left(\omega\right)$
 we denote by $\xi^{*}$ a representative point of $\xi$, and by
$G_{\xi}=\set{g\in G}{\omega\left(g\right)\xi^{*}=\xi^{*}}$ the stabilizer
of $\xi^{*}$. To each point $y\in\xi$ we associate a suitable element
$\gamma_{y}\in G$, such that $\xi^{*}$ is mapped to $y$ by the
permutation $\omega\left(\gamma_{y}\right)$ (of course, such a $\gamma_{y}$
is far from unique, any element of the coset $\gamma_{y}G_{\xi}$
would do the job).

\begin{lem}
\label{lem:Xhom}For a given homomorphism $\omega:G\rightarrow\Omega_{2}$,
the crossed homomorphisms $\lambda:G\rightarrow\Omega_{1}^{Y}$ \textup{}are
in one-to-one correspondence with pairs $\left(\Phi_{\xi},\varphi_{\xi}\right)$,
one for each orbit $\xi\in\mathcal{O}\left(\omega\right)$, where
$\Phi_{\xi}:G_{\xi}\rightarrow\Omega_{1}$ \textup{}is a homomorphism,
while $\varphi_{\xi}:\xi\rightarrow\Omega_{1}$ \textup{}is an arbitrary
map for which $\varphi_{\xi}\left(\xi^{*}\right)=\Phi_{\xi}\left(\gamma_{\xi^{*}}\right)$
.
\end{lem}
\begin{proof}
To begin with, let's note that any crossed homomorphism $\lambda:G\rightarrow\Omega_{1}^{Y}$
may be viewed as a map $\lambda:G\times Y\rightarrow\Omega_{1}$ that
satisfies \begin{equation}
\lambda\left(gh,y\right)=\lambda\left(g,hy\right)\lambda\left(h,y\right)\label{eq:xhom1}\end{equation}
for all $g,h\in G$ and $y\in Y$, where for simplicity $hy$ denotes
the image of $y$ under the permutation $\omega\left(h\right)\in\Omega_{2}$.
 By Eq.\eqref{eq:xhom1}, the map \begin{equation}
\begin{aligned}\Phi_{\xi}:G_{\xi} & \rightarrow\Omega_{1}\\
g & \mapsto\lambda\left(g,\FG{\xi}\right)\end{aligned}
\label{eq:Fixidef}\end{equation}
is a group homomorphism for each orbit $\xi\in\mathcal{O}\left(\omega\right)$.

Let's substitute $y=\FG{\xi}$ and $h=\gamma_{z}$ into Eq.\eqref{eq:xhom1}
to get\begin{equation}
\lambda\left(g,z\right)=\lambda\left(g\gamma_{z},\FG{\xi}\right)\lambda\left(\gamma_{z},\FG{\xi}\right)^{-1}\;\label{eq:xhom2}\end{equation}
for $z\in\xi$. By definition, $\left(g\gamma_{z}\right)\FG{\xi}=gz=\gamma_{gz}\FG{\xi}$,
in other words $\gamma_{gz}^{-1}g\gamma_{z}\in G_{\xi}$. Taking this
into account, Eq.\eqref{eq:xhom1} gives \begin{equation}
\lambda\left(g\gamma_{z},\FG{\xi}\right)=\lambda\left(\gamma_{gz},\FG{\xi}\right)\lambda\left(\gamma_{gz}^{-1}g\gamma_{z},\FG{\xi}\right)\:.\label{eq:xhom3}\end{equation}
If we introduce for each $\xi\in\mathcal{O}\left(\omega\right)$ the
maps \begin{equation}
\begin{aligned}\varphi_{\xi}:\xi & \rightarrow\Omega_{1}\\
y & \mapsto\lambda\left(\gamma_{y},\FG{\xi}\right)\:,\end{aligned}
\label{eq:varfidef}\end{equation}
 then Eqs.\eqref{eq:xhom2} and \eqref{eq:xhom3} lead to \begin{equation}
\lambda\left(g,y\right)=\varphi_{\xi}\left(gy\right)\Phi_{\xi}\left(\gamma_{gy}^{-1}g\gamma_{y}\right)\varphi_{\xi}\left(y\right)^{-1}\label{eq:xhom4}\end{equation}
for all $g\in G$ and $y\in\xi$. Note that, since $\gamma_{\xi^{*}}\in G_{\xi}$,
one has \begin{equation}
\varphi_{\xi}\left(\xi^{*}\right)=\lambda\left(\gamma_{\xi^{*}},\xi^{*}\right)=\Phi_{\xi}\left(\gamma_{\xi^{*}}\right)\:.\label{eq:xhom5}\end{equation}

Eq. \eqref{eq:xhom4} means that, given a choice of the coset representative
$\FG{\xi}$ and of the elements $\gamma_{y}$ mapping  $\FG{\xi}$
to $y\in\xi$, the pair $\left(\Phi_{\xi},\varphi_{\xi}\right)$ determines
completely the values $\lambda\left(g,y\right)$ for all points $y\in\xi$.
Since the orbits $\xi$ partition the set $Y$, it follows that the
collection of such pairs determines the crossed homomorphism $\lambda$.
\end{proof}
\begin{lem}
With the notations of Lemma \ref{lem:Xhom}, the orbits of the image
of the homomorphism \textup{$\phi:G\rightarrow\Omega_{1}\wr\Omega_{2}$}
are of the form $\left\langle \xi,\eta\right\rangle =\set{\left(\varphi_{\xi}\left(y\right)x,y\right)}{x\in\xi,\: y\in\eta}$,
where $\eta\in\mathcal{O}\left(\omega\right)$ and $\xi\in\mathcal{O}\left(\Phi_{\eta}\right)$\textup{.}
Moreover, the stabilizer of such an orbit is given by \textup{$G_{\left\langle \xi,\eta\right\rangle }=\set{g\in G_{\eta}}{\Phi_{\eta}\left(g\right)\xi^{*}=\xi^{*}}$.}
\end{lem}
\begin{proof}
Consider an arbitrary point $\left(x,y\right)\in X\times Y$: the
permutation $\phi\left(g\right)$ takes this to the point $\left(\lambda\left(g,y\right)x,gy\right)$:
in other words, the pair $\left(w,z\right)$ lies in the orbit of
$\left(x,y\right)$ if and only if there exists $g\in G$ such that
$z=gy$ and $w=\lambda\left(g,y\right)x$. By Eq.\eqref{eq:xhom4},
one has $\lambda\left(g,y\right)=\varphi_{\eta}\left(z\right)\Phi_{\eta}\left(\gamma_{z}^{-1}g\gamma_{y}\right)\varphi_{\eta}\left(y\right)^{-1}$
if $y$ lies in the orbit $\eta\in\mathcal{O}\left(\omega\right)$.
Now, for fixed $y$ and $z$ , the expression $\gamma_{z}^{-1}g\gamma_{y}$
runs through all elements of $G_{\eta}$. This shows that $\left(w,z\right)$
lies in the same orbit as $\left(x,y\right)$ if and only if $z$
lies in the $\omega\left(G\right)$-orbit $\eta$ of $y$, and at
the same time $\varphi_{\eta}\left(z\right)^{-1}w$ lies in the same
$\Phi_{\eta}\left(G_{\eta}\right)$-orbit as $\varphi_{\eta}\left(y\right)^{-1}x$,
which proves the claim about the structure of the orbits. The expression
for the stabilizer is obvious.
\end{proof}
Armed with the above results, we can now state and prove the fundamental
transitivity property of the orbifold transform.

\begin{thm}
Let $Z\in\FF{G}$ \textup{}be a class function of the finitely generated
group $G$, with values in the commutative $\mathbb{Q}$\textup{-}algebra
$R$, and let $\Omega_{1},\Omega_{2}$ be two permutation groups of
finite degree. Then \begin{equation}
\left(Z\wr\Omega_{1}\right)\wr\Omega_{2}=Z\wr\left(\Omega_{1}\wr\Omega_{2}\right)\;.\label{eq:transitivity}\end{equation}

\end{thm}
\begin{proof}
It is enough to prove that both sides of Eq.\eqref{eq:transitivity}
assign the same value to the subgroup $H=G$. Indeed, for a nontrivial
subgroup $H<G$, one may consider the restriction $Z_{H}\in\FF{H}$
of $Z$ to $\FE{H}$: since $Z_{H}$ agrees with $Z$ for all finite
index subgroups of $H$, it does so in particular for $H$ itself,
so Eq.\eqref{eq:transitivity} for $H$ holds if and only if it holds
for $Z_{H}$.

Let's first consider the rhs. of Eq.\eqref{eq:transitivity}, which
reads \begin{equation}
\frac{1}{\left|\Omega_{1}\wr\Omega_{2}\right|}\sum_{\phi:G\rightarrow\Omega_{1}\wr\Omega_{2}}\;\prod_{\zeta\in\mathcal{O}\left(\phi\right)}Z\left(G_{\zeta}\right)\:.\label{eq:transproof1}\end{equation}
By the previous arguments, we can associate to $\phi$ the pair $\left(\lambda,\omega\right)$,
where $\omega:G\rightarrow\Omega_{2}$ is a homomorphism, while $\lambda:G\times Y\rightarrow\Omega_{1}$
is a crossed homomorphism; moreover, $\lambda$ itself may be described
via pairs $\left(\Phi_{\eta},\varphi_{\eta}\right)$, one for each
orbit $\eta\in\mathcal{O}\left(\omega\right)$, according to Lemma
1. Taking this into account, as well as the structure of the orbits
of $\phi$ as described by Lemma 2, the expression in Eq.\eqref{eq:transproof1}
reads \begin{equation}
\frac{1}{\left|\Omega_{1}\right|^{\left|Y\right|}\left|\Omega_{2}\right|}\sum_{\omega:G\rightarrow\Omega_{2}}\:\prod_{\eta\in\mathcal{O}\left(\omega\right)}\:\sum_{\Phi_{\eta}:G_{\eta}\rightarrow\Omega_{1}}\:\sum_{\varphi_{\eta}}\:\prod_{\xi\in\mathcal{O}\left(\Phi_{\eta}\right)}Z\left(G_{\left\langle \xi,\eta\right\rangle }\right)\:.\label{eq:transproof2}\end{equation}
This may be rearranged into the more suggestive form \begin{equation}
\frac{1}{\left|\Omega_{2}\right|}\sum_{\omega:G\rightarrow\Omega_{2}}\:\prod_{\eta\in\mathcal{O}\left(\omega\right)}\left\{ \sum_{\varphi_{\eta}}\frac{1}{\left|\Omega_{1}\right|^{\left|\eta\right|}}\sum_{\Phi_{\eta}:G_{\eta}\rightarrow\Omega_{1}}\;\prod_{\xi\in\mathcal{O}\left(\Phi_{\eta}\right)}Z\left(G_{\left\langle \xi,\eta\right\rangle }\right)\right\} \:.\label{eq:transproof3}\end{equation}
Nothing in this expression depends explicitly on the maps $\varphi_{\eta}$,
so their contribution is simply to introduce a multiplicative factor
$\left|\Omega_{1}\right|^{\left|\eta\right|-1}$ equal to their number.
Taking into account the structure of $G_{\left\langle \xi,\eta\right\rangle }$
as described by Lemma 2, we recognize that the resulting expression
is nothing but \begin{equation}
\frac{1}{\left|\Omega_{2}\right|}\sum_{\omega:G\rightarrow\Omega_{2}}\prod_{\eta\in\mathcal{O}\left(\omega\right)}\left(Z\wr\Omega_{1}\right)\left(G_{\eta}\right)\:,\label{eq:transproof4}\end{equation}
which is just the lhs. of Eq.\eqref{eq:transitivity}.
\end{proof}
We note that it is this transitivity property Eq.\eqref{eq:transitivity}
that motivates the wreath product notation $Z\wr\Omega$ for the orbifold
transform (besides the connection with permutation orbifolds).

\section{The exponential identity and symmetric products}

In Section 2, when discussing examples of the orbifold transform,
we have fixed the finitely generated group $G$, and let both the
class function $Z$ and the permutation group $\Omega$ vary freely.
Another possible approach is to fix the permutation group $\Omega$,
while leaving $G$ and $Z$ arbitrary: a most important case is when
$\Omega=S_{n}$, the symmetric group of degree $n$, which is termed
a symmetric product in the theory of permutation orbifolds. As it
turns out, one has very good control over the orbifold transform in
this case, thanks to a general combinatorial identity \cite{symprod}
that we are going to discuss.

First, let's fix some notation. For  a positive integer $n$ and an
arbitrary class function $Z\in\FF{G}$ of the finitely generated group
$G$, we define $Z_{n}=Z\wr S_{n}$, and we set $Z_{0}$ to the constant
class function equal to $1$: we call $Z_{n}$ the $n$-th symmetric
product of $Z$.

\begin{thm}
The following formal identity holds

\begin{equation}
\sum_{n=0}^{\infty}Z_{n}\left(G\right)=\exp\left(\sum_{H\in\FE{G}}\frac{Z\left(H\right)}{\left[G:H\right]}\right)\;,\label{eq:expoid}\end{equation}
where the exponential on the rhs. stands for its infinite power series.
\end{thm}
\begin{proof}
Let's write out the lhs. of Eq.\eqref{eq:expoid}: it reads \begin{equation}
1+\sum_{n=1}^{\infty}\frac{1}{n!}\sum_{\phi:G\rightarrow S_{n}}\:\prod_{\xi\in\mathcal{O}\left(\phi\right)}Z\left(G_{\xi}\right)\:.\label{eq:expproof1}\end{equation}
Let's first consider the sum over homomorphisms $\phi:G\rightarrow S_{n}$
for a given $n$: each such homomorphism is just a permutation action
of $G$ of degree $n$, and $\mathcal{O}\left(\phi\right)$ is the
set of orbits of this action, while the $G_{\xi}$ are the point stabilizers.
The point is that the term corresponding to a given $\phi$ does only
depend on the equivalence class of this action, since there is a one-to-one
correspondence between the orbits of equivalent actions, and the corresponding
stabilizers are conjugate subgroups. This means that we can rewrite
this sum over all permutation actions of degree $n$ as a sum over
equivalence classes of actions, provided we take into account the
cardinality of each equivalence class.

Any equivalence class of permutation actions may be decomposed into
a sum of transitive classes, corresponding to the different orbits
of the action: \begin{equation}
\left[\phi\right]=\oplus_{i}n_{i}\tau_{i}\:,\label{eq:phidecomp}\end{equation}
where the $\tau_{i}$ denote the different transitives, and $n_{i}$
is the multiplicity of the $i$-th transitive in $\left[\phi\right]$.
For each given degree there are only finitely many transitives of
that degree: indeed, transitives correspond to conjugacy classes of
subgroups of finite index, the degree of the transitive being equal
to the index of the subgroup.

Let $G_{i}$ denote the stabilizer of the $i$-th transitive $\tau_{i}$:
in other words, $\tau_{i}$ is the equivalence class of the action
of $G$ on the cosets of the subgroup $G_{i}$. Then, for each permutation
action in the equivalence class $\left[\phi\right]$, one has \begin{equation}
\prod_{\xi\in\mathcal{O}\left(\phi\right)}Z\left(G_{\xi}\right)=\prod_{i}Z\left(G_{i}\right)^{n_{i}}\:.\label{eq:expproof2}\end{equation}
Denoting by $\FD=\left[N_{G}\left(G_{i}\right):G_{i}\right]$ the
index of $G_{i}$ in its normalizer, the cardinality of the equivalence
class $\left[\phi\right]$ is given by \begin{equation}
\#\left[\phi\right]=\frac{n!}{\prod_{i}n_{i}!\FD^{n_{i}}}\:,\label{eq:phiclassnumber}\end{equation}
where $n$ is the degree of $\left[\phi\right]$ (see Appendix A for
a detailed proof). Taking all this into account,  the rhs. of Eq.\eqref{eq:expoid}
reads \begin{equation}
1+\sum_{\left[\phi\right]}\prod_{i}\frac{1}{n_{i}!\FD^{n_{i}}}Z\left(G_{i}\right)^{n_{i}}\:,\label{eq:expproof3}\end{equation}
where the summation runs over all equivalence classes $\left[\phi\right]$.
But the multiplicities $n_{i}$ may take on arbitrary nonnegative
values, which leads to \begin{equation}
\sum_{n=0}^{\infty}Z_{n}\left(G\right)=\prod_{i}\left(\sum_{n_{i}=0}^{\infty}\frac{Z\left(G_{i}\right)^{n_{i}}}{n_{i}!\FD^{n_{i}}}\right)=\prod_{i}\exp\left(\frac{Z\left(G_{i}\right)}{\FD}\right)\:,\label{eq:expproof4}\end{equation}
where the product runs over all transitives of finite degree, or --
what is the same -- over all conjugacy classes of finite index subgroups.
Since the number of different conjugates of $G_{i}$ equals \begin{equation}
\left[G:N_{G}\left(G_{i}\right)\right]=\frac{\left[G:G_{i}\right]}{\FD}\:,\label{eq:expproof5}\end{equation}
and a product of exponentials is the exponential of the sum of the
exponents, we get the assertion of the theorem.
\end{proof}
The exponential identity is formal in the sense that convergence of
the infinite sums is not guaranteed on either side. Of course, for
suitable choice of the class function $Z$ one obtains an equality
of convergent series.

Let us rewrite the exponential identity Eq.\eqref{eq:expoid} into
a slightly different form. To this end, we introduce the notation

\begin{equation}
Z^{\left[n\right]}\left(G\right)=\sum_{\left[G:H\right]=n}Z\left(H\right)\:,\label{eq:genHecke}\end{equation}
for a positive integer $n$, where the summation on the rhs. extends
over all subgroups $H<G$ of index $n$. Moreover, we adjoin a formal
variable $p$ to the ring $R$, and consider the class function (where
$R\left\{ p\right\} $ stands for the ring of formal power series
in $p$ with coefficients from $R$) \begin{eqnarray*}
\hat{Z}:\FE{G} & \rightarrow & R\left\{ p\right\} \\
H & \mapsto & p^{\left[G:H\right]}Z\left(H\right)\:.\end{eqnarray*}
Applying the exponential identity Eq.\eqref{eq:expoid} to the class
function $\hat{Z}$, one gets the following: \begin{equation}
\sum_{n=0}^{\infty}p^{n}Z_{n}\left(G\right)=\exp\left(\sum_{n=1}^{\infty}\frac{p^{n}Z^{\left[n\right]}\left(G\right)}{n}\right)\:.\label{eq:expoid2}\end{equation}

Now, denoting by $P_{n}$ the cycle indicator of the symmetric group
$S_{n}$ (a Schur-polynomial), one has the following well-known identity
\cite{stanley}: \begin{equation}
1+\sum_{n=1}^{\infty}p^{n}P_{n}\left(t_{1},\ldots,t_{n}\right)=\exp\left(\sum_{n=1}^{\infty}\frac{p^{n}t_{n}}{n}\right)\:,\label{eq:Snexpo}\end{equation}
which is actually a special case of Eq.\eqref{eq:expoid} (for $G=\mathbb{Z}$,
$R=\mathbb{Q}\left[t_{1},t_{2},\ldots\right]$ and $Z\left(k\mathbb{Z}\right)=t_{k}$,
cf. Section 2). Comparing Eqs.\eqref{eq:expoid2} and \eqref{eq:Snexpo},
and equating the coefficients of equal powers of $p$, one arrives
at\begin{equation}
Z_{n}\left(G\right)=P_{n}\left(Z^{\left[1\right]}\left(G\right),\ldots,Z^{\left[n\right]}\left(G\right)\right)\:,\label{eq:symprod}\end{equation}
providing an elegant closed formula for the symmetric products of
$Z$.

\section{The transform and permutation orbifolds}

Consider a system that is made up of several identical, non-interacting
subsystems, each of which may be described by a  Conformal Field Theory.
The dynamics of the whole system is again governed by a CFT, but what
is more important, any permutation of the subsystems is a symmetry
of this CFT, since the subsystems are indistinguishable: consequently,
one may orbifoldize with respect to any group $\Omega$ of permutations
of the subsystems. The resulting CFT, which is completely determined
by the twist group $\Omega$ and the CFT $\mathcal{C}$ describing
the dynamics of the individual subsystems, is the permutation orbifold
$\mathcal{C}\wr\Omega$.

Many important aspects %
 of permutation orbifolds are pretty well understood: one knows how
to classify their primary fields, one has elegant closed expressions
for the genus one characters of these primaries, the modular transformations
of the characters, the fusion rules, the partition functions, etc
\cite{permorbi1,permorbi2,uniformization}. This is where the orbifold
transform enters the picture: we'll illustrate this point on the example
of partition functions.

Among other things, a CFT assigns a number to each conformal equivalence
class of two-dimensional metrics: this is what is called the partition
function (more precisely, the generalized partition function, the
usual partition function is obtained by restricting attention to metrics
defined on tori). It is well known that for orientable surfaces, conformal
equivalence classes of metrics  are in one-to-one correspondence with
equivalence classes of complex structures, i.e. the partition function
may be viewed as a function on the moduli space of Riemann surfaces.
For simplicity, we'll restrict our attention to closed compact %
 surfaces.

By the uniformization theorem \cite{Beardon,farkas-kra}, every Riemann
surface $\mathcal{S}$ may be obtained as a quotient \begin{equation}
\hat{\mathcal{S}}/G_{\mathcal{S}}\:,\label{eq:uniformization}\end{equation}
where $\hat{\mathcal{S}}$ is a simply connected Riemann surface (the
universal cover of $\mathcal{S}$), while $G_{\mathcal{S}}$, the
uniformizing group of $\mathcal{S}$, is a discrete group of automorphisms
of $\hat{\mathcal{S}}$, isomorphic to the fundamental group $\pi_{1}\left(\mathcal{S}\right)$
of $\mathcal{S}$. Actually, the uniformizing group is only determined
up to conjugacy in $\aut{\hat{\mathcal{S}}}$, i.e. conjugate subgroups
uniformize the same surface.

By the Riemann mapping theorem \cite{Beardon,farkas-kra}, there are
just three inequivalent simply connected Riemann surfaces: the Riemann
sphere $\mathbb{CP}^{1}=\mathbb{C}\cup\left\{ \infty\right\} $, the
complex plane $\mathbb{C}$  and the upper half-plane $\FB=\set{z}{\mathrm{Im}\!\ z>0}$.
Among compact closed surfaces, the only one whose universal cover
is $\mathbb{CP}^{1}$ is the Riemann sphere itself, and the uniformizing
group is trivial; $\mathbb{C}$ is the universal cover of tori, and
in this case the uniformizing group is a group of complex translations
isomorphic to $\FA$; finally, $\FB$ is the universal cover of all
other compact surfaces, and these are uniformized by  Fuchsian groups.
In terms of the genus $g$ of the surface $\mathcal{S}$, the above
cases correspond respectively to $g=0$, $g=1$ and $g>1$.

The theory of covering surfaces \cite{Forster,Beardon} tells us that
any finite index subgroup $H<G_{\mathcal{S}}$ is again the uniformizing
group of some compact Riemann surface $\hat{\mathcal{S}}/H$, which
is a finite sheeted cover of $\mathcal{S}$. Since uniformizing groups
are finitely generated (being isomorphic to the fundamental group
of a surface of finite genus), given a CFT and compact Riemann surface
$\mathcal{S}$, we can define a class function $\mathcal{Z_{\mathcal{S}}}\in\mathscr{C}\!\left(G_{\mathcal{S}},\mathbb{C}\right)$
by assigning to each subgroup $H\in\FE{G_{\mathcal{S}}}$ the value
of the partition function on the surface $\hat{\mathcal{S}}/H$: the
facts from complex analysis listed above ensure that this is well
defined.

Going back to permutation orbifolds, the CFT $\mathcal{C}$ that describes
the dynamics of the individual subsystems gives rise, according to
the preceding discussion, to a class function $\mathcal{Z_{\mathcal{S}}}\in\mathscr{C}\!\left(G_{\mathcal{S}},\mathbb{C}\right)$
for any compact surface $\mathcal{S}$. The permutation orbifold $\mathcal{C}\wr\Omega$
leads to another class function $\mathcal{Z_{\mathcal{S}}^{\mathnormal{\mathrm{\Omega}}}}\in\mathscr{C}\!\left(G_{\mathcal{S}},\mathbb{C}\right)$,
and the question is whether these two class functions are related
or not. It follows from the results of \cite{uniformization}, that
for any compact surface $\mathcal{S}$ the class function $\mathcal{Z_{\mathcal{S}}^{\mathnormal{\mathrm{\Omega}}}}$
is the orbifold transform of $\mathcal{Z_{\mathcal{S}}}$: \begin{equation}
\mathcal{Z_{\mathcal{S}}^{\mathnormal{\mathrm{\Omega}}}}=\mathcal{Z_{\mathcal{S}}}\wr\Omega\:.\label{eq:genpartfun}\end{equation}
This is the basic connection with permutation orbifolds. It should
be stressed that this connection is not confined to partition functions,
many other characteristic quantities of permutation orbifolds may
be expressed as suitable orbifold transforms, e.g. the number of primary
fields, the traces of mapping class group transformations, etc \cite{permorbi2}.

Of course, the result that the partition function of the permutation
orbifold is the orbifold transform of the partition function of a
single subsystem is far from being trivial. It is based on the  physical
picture that the dynamics of the orbifold on a given world sheet may
be interpreted as the dynamics of a subsystem on suitable covers of
the world sheet \cite{uniformization}. In this respect, the transitivity
property Eq.\eqref{eq:transitivity} plays a decisive role: indeed,
the corresponding property of permutation orbifolds is an immediate
consequence of their definition \cite{permorbi1}, which has to manifest
itself in any expression relating the quantities of the orbifold with
those of the subsystems. Actually, Eq.\eqref{eq:transitivity} would
have been difficult to guess were it not for the connection with permutation
orbifolds.

To conclude this section, let's illustrate the above results on the
simplest nontrivial case, the genus one partition function $Z$. It
is defined on the moduli space of two dimensional tori, which is just
the quotient of the upper half-plane $\FB$ by the classical modular
group $\FC$. It is usually written in terms of the modular parameter
$\tau\in\FB$ of the corresponding torus, in terms of which it satisfies
the functional equation (modular invariance) \begin{equation}
Z\left(\frac{a\tau+b}{c\tau+d}\right)=Z\left(\tau\right)\:,\label{eq:modinv}\end{equation}
for all $\sm{a}{b}{c}{d}\in\FC$.

The universal cover of the torus with modular parameter $\tau$ is
the complex plane, and its uniformizing group may be taken to be the
group $G_{\tau}$ generated by the two translations \begin{equation}
\begin{aligned}a:z & \mapsto z+1\\
b:z & \mapsto z+\tau\:.\end{aligned}
\label{eq:torusgens}\end{equation}
Clearly,  $G_{\tau}$ is isomorphic to $\FA$, which is the last example
discussed in Section 2. A finite index subgroup of $G_{\tau}$ corresponds
to a matrix \begin{equation}
H=\left(\begin{array}{cc}
\mu & \kappa\\
0 & \lambda\end{array}\right)\label{eq:torusHNF}\end{equation}
in HNF, and is generated by \begin{equation}
\begin{aligned}a^{\lambda}:z & \mapsto z+\lambda\\
a^{\kappa}b^{\mu}:z & \mapsto z+\mu\tau+\kappa\:.\end{aligned}
\label{eq:subgroupgens}\end{equation}
But this subgroup is conjugate in $\aut{\mathbb{C}}$ to the group
$G_{\tau{\scriptstyle \left(H\right)}}$, where \begin{equation}
\tau\left(H\right)=\frac{\mu\tau+\kappa}{\lambda}\:,\label{eq:subgrouptau}\end{equation}
i.e. it uniformizes a torus with modular parameter $\tau\left(H\right)$.

Taking this into account, as well as Eqs.\eqref{eq:torustransform}
and \eqref{eq:genpartfun}, we arrive at the following well-known
expression for the torus partition function of permutation orbifolds
\cite{permorbi1}: \begin{equation}
Z^{\Omega}\left(\tau\right)=\frac{1}{\left|\Omega\right|}\sum_{{{x,y\in\Omega\atop xy=yx}}}\:\prod_{\xi\in\mathcal{O}\left(x,y\right)}Z\left(\frac{\mu_{\xi}\tau+\kappa_{\xi}}{\lambda_{\xi}}\right)\:,\label{eq:toruspartfun}\end{equation}
where the integers $\lambda_{\xi},\mu_{\xi}$ and $\kappa_{\xi}$
are numerical characteristics of the orbit $\xi\in\mathcal{O}\left(x,y\right)$,
whose meaning is described just before Eq.\eqref{eq:torustransform}.

\section{Summary and outlook}

The aim of this note was to introduce the notion of the orbifold transform,
to illustrate it on some simple examples, and to discuss its most
basic properties and some of its applications. The two major results,
namely the transitivity property Eq.\eqref{eq:transitivity} of the
transform and the exponential identity Eq.\eqref{eq:expoid} already
justify amply the consideration of this construct, and the connection
with the theory of permutation orbifolds gives even more evidence
of its importance.

One should mention that it is possible to generalize the transform
by including suitable 'cohomological twists': these arise naturally
in the CFT context, where they run under the name of discrete torsion.
The resulting theory is pretty similar to the one described here,
but we refrained from its exposition, since it would need a thorough
treatment of the cohomology of wreath products, which would be beyond
the scope of this note. While most results go over to this general
case, some of them get modified, e.g. the exponential formula Eq.\eqref{eq:expoid}:
for those interested in this issue we mention \cite{symprod}, where
the torus partition function of symmetric products in the presence
of nontrivial discrete torsion is discussed.

\appendix

\section{~}

This appendix is devoted to the proof of Eq.\eqref{eq:phiclassnumber},
giving the number of different permutation actions in the same equivalence
class. First, let's fix the notation. Let $\tau_{i}$ denote the different
equivalence classes of transitive actions of $G$; let $G_{i}$ denote
the stabilizer of some point of $\tau_{i}$ (the stabilizers of different
points form a conjugacy class of subgroups), so that $\tau_{i}$ is
equivalent to the action of $G$ on the cosets of $G_{i}$; finally,
let $L_{i}$ denote the factor group $N_{G}\left(G_{i}\right)/G_{i}$.

\begin{lem}
Let $\phi:G\rightarrow S_{n}$ \textup{}be a permutation action. The
number of different permutation actions equivalent to $\phi$ \textup{}equals
the index of the centralizer $C\!\left[\phi\right]$ \textup{}of the
image of $\phi$ \textup{}in \textup{$S_{n}$:} \begin{equation}
\#\left[\phi\right]=\left[S_{n}:C\!\left[\phi\right]\right]\:.\label{eq:appfi}\end{equation}

\end{lem}
\begin{proof}
The equivalence class $\left[\phi\right]$ consists of the conjugates
\begin{equation}
\begin{aligned}\phi^{\alpha}:G\rightarrow & S_{n}\\
g\mapsto & \alpha^{-1}\phi\left(g\right)\alpha\end{aligned}
\label{eq:appfi2}\end{equation}
of $\phi$, where $\alpha\in S_{n}$. Clearly, the assignment $\phi\mapsto\phi^{\alpha}$
defines a permutation action of $S_{n}$ on the set of actions of
$G$ of degree $n$, and $\left[\phi\right]$ is the orbit of $\phi$
under this action. The length of the orbit equals the index of the
stabilizer of $\phi$, but this stabilizer is nothing but the centralizer
$C\!\left[\phi\right]$ of the image of $\phi$.
\end{proof}
\begin{lem}
If the permutation action $\phi$ \textup{}belongs to the equivalence
class $\oplus_{i}n_{i}\tau_{i}$, then the centralizer of its image
is isomorphic to \begin{equation}
\times_{i}\left(L_{i}\wr S_{n_{i}}\right)\:,\label{eq:centstructure}\end{equation}
\textup{i}n particular Eq.\eqref{eq:phiclassnumber} holds.
\end{lem}
\begin{proof}
By the definition of equivalence of permutation actions, one has a
direct product decomposition \begin{equation}
C\left[\oplus_{i}n_{i}\tau_{i}\right]=\times_{i}C\left[n_{i}\tau_{i}\right]\:.\label{eq:centproof1}\end{equation}
Now,  the centralizer $C\left[n_{i}\tau_{i}\right]$ is generated
by two sorts of permutations: those that permute the transitive constituents
(i.e. the orbits) en block, without permuting the points inside an
orbit, which form a group isomorphic to $S_{n_{i}}$; and those that
do leave each orbit setwise fixed, but permute the points of the orbits
while still commuting with the action. These two sort of permutations
generate the wreath product \begin{equation}
C\left[n_{i}\tau_{i}\right]=C\left[\tau_{i}\right]\wr S_{n_{i}}\:.\label{eq:centproof2}\end{equation}
We have reduced the problem to that of determining $C\left[\tau_{i}\right]$
for a transitive action. But such a transitive action is equivalent
to the action of $G$ on the cosets of $G_{i}$, so we are looking
for permutations $\alpha\in S_{G/G_{i}}$ such that \begin{equation}
\alpha\left(gxG_{i}\right)=g\alpha\left(xG_{i}\right)\label{eq:centcondition}\end{equation}
 holds for all $g,x\in G$ , which is just the condition $\alpha\in C\left[\tau_{i}\right]$.
Eq.\eqref{eq:centcondition} with $x=1$ gives $\alpha\left(gG_{i}\right)=g\alpha\left(G_{i}\right)$,
i.e. the permutation $\alpha$ is completely determined by the image
$\alpha\left(G_{i}\right)$ of the trivial coset. Since this image
is itself a coset of $G_{i}$, there exists some $a\in G$ such that
$\alpha\left(G_{i}\right)=aG_{i}$. Now, Eq.\eqref{eq:centcondition}
with $g\in G_{i}$ yields \begin{equation}
aG_{i}=\alpha\left(G_{i}\right)=g\alpha\left(G_{i}\right)=gaG_{i}\:,\label{eq:centproof3}\end{equation}
in other words $a$ should belong to the normalizer of $G_{i}$ in
$G$ in order to get a permutation of the cosets. This means that
\begin{equation}
C\left[\tau_{i}\right]\cong N_{G}\left(G_{i}\right)/G_{i}\:.\label{eq:centproof4}\end{equation}
All-in-all, we arrive at Eq.\eqref{eq:centstructure}. Combining this
with Eq.\eqref{eq:appfi}, we get Eq.\eqref{eq:phiclassnumber}.
\end{proof}
\bibliographystyle{plain}

\end{document}